\documentclass[12pt]{article}
\usepackage[all]{xy}
\usepackage{wrapfig}

\usepackage{epic}
\usepackage{amsmath}[1996/11/01]
\usepackage{amssymb,amsthm,amsfonts,latexsym,epsfig,graphics}
 \def\beql#1#2\eeql{\begin{equation}\label{#1}#2\end{equation}}

\advance \textheight by -1 \topmargin
\advance \textheight by -1 \headheight
\advance \textheight by -1 \headsep
\advance \textheight by -2 \footskip
\vsize = \textheight
\hsize = \textwidth
\DeclareMathOperator{\rad}{rad}
\DeclareMathOperator{\Gal}{Gal}
\DeclareMathOperator{\Deltaa}{{\mathfrak D}}
\DeclareMathOperator{\Rho}{{\mathfrak R}}

\DeclareMathOperator{\Irr}{Irr}

\DeclareMathOperator{\Char}{char}
\DeclareMathOperator{\disc}{disc}

\DeclareMathOperator{\diag}{diag}

\theoremstyle{plain}
\newtheorem{theorem}{Theorem}
\newtheorem{lemma}[theorem]{Lemma}

\newtheorem{conjecture}[theorem]{Conjecture}

\newtheorem{proposition}[theorem]{Proposition}
\newtheorem{corollary}[theorem]{Corollary}
\newtheorem{definition}[theorem]{Definition}
\theoremstyle{remark}

\newtheorem{remark}[theorem]{Remark}

\numberwithin{theorem}{section}

\newcommand{\Z}{{\mathbb{Z}}}
\newcommand{\Q}{{\mathbb{Q}}}

\newcommand{\N}{{\mathbb{N}}}

\renewcommand{\em}{\sf}

\title{On orthogonal discriminants of characters}
\author{Gabriele Nebe\footnote{nebe@math.rwth-aachen.de}}  
\date{Lehrstuhl f\"ur Algebra und Zahlentheorie, 
RWTH Aachen University, Germany}
 
\begin{document}
\maketitle

 {\sc Abstract.} 
An ordinary character $\chi $ of a finite group is called 
 orthogonally stable, if all 
non-degenerate 
 invariant quadratic forms on any module affording the character $\chi $ 
 have the same discriminant. This is the 
 orthogonal discriminant, $\disc(\chi )$, of $\chi $, a square class of the 
 character field. 
 Based on experimental evidence we conjecture that the orthogonal 
 discriminant is always an odd square class in the sense of 
 Definition \ref{oddsc}.
This note  proves this conjecture for finite solvable groups. 
For $p$-group there is an explicit formula for 
 $\disc(\chi )$ that reads 
 $\disc(\chi ) = (-p)^{\chi(1)/2}$ if $p\equiv 3 \pmod{4}$ and 
 $\disc (\chi ) = (-1)^{\chi(1)/2}$ for $p=2$.
     \\
  {\sc keywords:}  orthogonal representations of finite groups; 
  character fields; discriminant fields; orthogonal discriminants.

\section{Introduction} 

In the past year Richard Parker and the author started a long term
project studying quadratic forms invariant under the finite simple groups 
in the Atlas of Finite Groups \cite{ATLAS}. 
We work over number fields  and 
finite fields and use decomposition matrices
to compute 
 discriminants of invariant quadratic forms (see for instance \cite[Section 6]{OrthogonalStability} for explicit examples).

 If $\Char (K) \neq 2$ then 
 any quadratic form $Q:V\to K$ on a $K$-space $V$ is determined 
by the corresponding bilinear form 
$$B:V\times V \to K, B(x,y) = Q(x+y) - Q(x) - Q(y) .$$ 
 The quadratic form is called non-degenerate, if the radical of $B$ is $\{0 \}$
 where
$$\rad(B) = V^{\perp} = \{ x\in V \mid B(x,y) = 0 \mbox{ for all } y\in V \} .$$
The discriminant of a non-degenerate quadratic form $Q$ is 
$$\disc((V,Q)) = \disc(Q) = \disc(B) = 
 (-1)^{{n}\choose{2}} \det (B) \in K^{\times}/(K^{\times})^2$$  where $n = \dim (V)$ and 
$\det (B)$ is the square class of the  determinant of a Gram matrix of $B$ with respect
to a chosen basis of $V$. 
The discriminant is a strong invariant of the isometry class of 
a non-degenerate quadratic form. 
For finite fields of odd characteristic discriminant and dimension 
completely determine the isometry class of a quadratic space
 (see \cite[Section IV]{Kneser}). 

The present short note is mainly concerned with number fields $K$.
So to simplify notation we assume that $K$ is a finite extension of 
the rational field $\Q $. Let $V$ be a $KG$-module. 
Then $V$ is called {\em orthogonal} if 
$G$ fixes a non-degenerate quadratic form on $V$.
 An orthogonal $KG$-module $V$ is called {\em orthogonally stable}
if all non-degenerate $G$-invariant quadratic forms  have 
the same discriminant. This square class of $K$ is then called
the {\em orthogonal discriminant} of the $KG$-module $V$. 

The properties of being orthogonal and orthogonally stable 
can be read off from the character of $V$: 
Let $\Irr(G)$  denote the set of absolutely 
irreducible ordinary characters of $G$. 
For $\psi \in \Irr(G)$ the Frobenius Schur indicator of $\psi $ 
 (for short {\em indicator} of $\psi$) takes the values 
 $o$, $+$, and $-$. 
 The indicator of $\psi $ 
is $o$, if $\psi $ takes non-real values. 
If $\psi $ is the character of a real representation then
the indicator of $\psi $ is $+$ and the indicator is
 $-$ if $\psi $ is real valued but there is no real 
representation affording the character $\psi $. 
Any character $\chi $ of $G$ is a unique sum of absolute irreducible 
 characters.  
 Recall that the character field of $\chi $ is the 
 number field generated by the character values, 
 $\Q(\chi ) := \Q (\chi (g) : g \in G ) $. 
Now $\chi (g)$ is the trace of the endomorphism of $V$ defined by $g$, 
so for 
 any $KG$-module $V$ affording the character $\chi $ the
 character field is a subfield of $K$.  
 If the Schur index of $\chi $ is 
 bigger than 1  then there is no $\Q(\chi ) G$-module with character $\chi $.

 \begin{proposition} (see \cite{OrthogonalStability})
	 A character 
	 $$\chi := \sum _{\psi \in \Irr(G) } n_{\psi } \psi  $$
	 is {\em orthogonal} if $\chi $ is real valued 
 and $n_{\psi }$ is even for all 
 $\psi \in \Irr(G)$ with  indicator $-$. 
	 \\
	 If $\chi $ is orthogonal then it is {\em orthogonally stable}
	 if $\psi (1) $ is even for all $\psi \in \Irr (G)$ of 
	 indicator $+$ for which $n_{\psi } > 0$.
 \end{proposition}

 The papers \cite{orthdet} and \cite{OrthogonalStability} 
 show that the orthogonal discriminant can be defined 
 as a square class of the character field: 
 Given an orthogonally stable 
 character $\chi $ there is a unique square class $\disc(\chi ) \in 
 \Q(\chi)^{\times}/(\Q(\chi)^{\times})^2 $ 
 such that for any  $KG$-module $V$ affording
 $\chi $ all non-degenerate $G$-invariant 
 quadratic forms on $V$  have discriminant
 $ \disc(\chi ) (K^{\times})^2$. 

 \begin{definition}
	 $\disc(\chi ) \in \Q(\chi)^{\times}/(\Q(\chi)^{\times})^2 $ 
	 is called the {\em orthogonal discriminant} 
	 of the orthogonally stable character $\chi $.
 \end{definition}

Richard Parker, Jon Thackray, Thomas Breuer  and
the author computed thousands orthogonal 
discriminants of orthogonally stable 
characters of finite simple groups. 
For rational characters the orthogonal discriminants are
represented by a unique square free integer $d$. 
We never found an even number $d$ which let Richard Parker 
to formulate the following conjecture. 

\begin{conjecture}
Let $\chi $ be an orthogonally stable rational character.
Then the orthogonal discriminant of $\chi $ is represented 
by an odd square free integer. 
\end{conjecture}

From the classification
of all maximal finite rational matrix groups up to dimension 31 
 (\cite{NeP}, \cite{dim24}, \cite{dim2531}) one can check that the 
conjecture is true if $\chi (1) \leq 30$. 
The aim of the present note is to prove a more general result
(Theorem \ref{main})
for all orthogonally stable characters of finite solvable groups.
To formulate the result for arbitrary character fields 
we need to have an appropriate notion of an odd square class.
So let $K$ be a number field and $\nu : K^{\times} \to \Z $ a 
surjective discrete valuation. 
Then  for $\delta \in K^{\times}$
$$\nu (\delta (K^{\times})^2) = \nu (\delta ) + 2\Z $$ 
so any square class in $K^{\times}$ has either even or odd valuation. 
We call $\nu $ a {\em dyadic} valuation, if $\nu (2) > 0$. 
\begin{definition}\label{oddsc}
	A square class  $\delta (K^{\times})^2$ is called odd, if $\nu (\delta (K^{\times})^2) $ 
	is even for all dyadic valuations $\nu : K^{\times} \to \Z $.
\end{definition}

Then the main result of this short note is the
following theorem:

\begin{theorem}\label{main}
Let $G$ be a finite solvable group and $\chi $ be 
an orthogonally stable ordinary character of $G$. 
Then $\disc(\chi) $ is  odd.
\end{theorem}

We also compute the orthogonal discriminants for
characters of $p$-groups in Theorem \ref{mainP} and \ref{main2} 
below. In particular we get that 
  the orthogonal discriminant of any orthogonally stable character
  $\chi $ of a $p$-group is $(-p)^{\chi (1)/2 } $ 
 if $p\equiv 3 \pmod{4} $ 
 and $(-1)^{\chi(1)/2} $ if $p=2$. 
 For primes $p\equiv 1 \pmod{4}$
 Theorem \ref{mainP} gives a similarly explicit result.

 \section{Preliminaries} 

 \subsection{The main character theoretic tool} 

It is well known that the character field $\Q(\chi )$ of any 
ordinary character $\chi $ of some finite group $G$ is an abelian
extension of $\Q $ and hence contained in some 
cyclotomic field $\Q[\zeta _f]$. The minimal such $f$ is 
called the conductor of $\chi $. 

\begin{theorem}(\cite[Theorem A1]{NavarroTiep}) \label{odddeg}
Let $\psi $ be an absolutely irreducible 
ordinary character of some finite group $G$ 
such that $\psi(1) $ is odd. 
Write the conductor of $\psi $ as $2^a m$ with $m\in \N$ odd. 
Then $\Q[\zeta _{2^a}] \subseteq \Q(\psi )$.
\end{theorem}

For our purposes a weaker version of Theorem \ref{odddeg} suffices:

\begin{theorem} \label{odddegcor} (\cite[Theorem C]{ILNT}) 
	Let $\psi $ be an absolutely irreducible
ordinary character of some finite group $G$
such that $\psi(1) $ is odd.
Then either the conductor of $\psi $ is odd or 
	$\Q[\sqrt{-1}] \subseteq \Q(\psi )$.
\end{theorem}

\begin{corollary} \label{realodd}
The conductor of a real, absolutely irreducible, odd degree character
 is odd.
\end{corollary} 

In Section \ref{pgrps} orthogonal discriminants of $p$-groups are 
determined.
The proof uses the fact that  for odd primes $p$ the
character fields of absolutely irreducible characters of $p$-groups are  
cyclotomic number fields: 

\begin{theorem} (see for instance \cite[Theorem 2.3]{NavarroTiep}) \label{pgrpchar}
	Let $p$ be an odd prime, $G$ a $p$-group and $1\neq \chi \in \Irr(G)$ 
	be a non-trivial absolutely irreducible character of $G$.
	Then $\Q(\chi ) = \Q [\zeta _{p^f}]$ for some $f\in \N$. 
\end{theorem} 

\subsection{The discriminant field}

Let $K$ be a number field and $\delta \in K^{\times}$. 
Then the square class $\delta (K^{\times})^2 $ defines a unique field extension 
$K[\sqrt{\delta }]$ of $K$. Note that $K[\sqrt{\delta }] = K$ if and only if 
$\delta $ is a square in $K$.

Certain results for quadratic forms over number fields have a more natural formulation
if one replaces the discriminant by the discriminant field:

\begin{definition} 
	 Let $\chi $ be an orthogonally stable character with 
		character field $K=\Q(\chi)$ and orthogonal discriminant
		 $\disc(\chi ) = \delta (K^{\times})^2$.
		Then 
		the {\em discriminant field} of $\chi $ is
                    $\Deltaa (\chi ) := K[\sqrt{\delta }] $.
\end{definition}

Whereas character fields are always abelian number fields and in particular Galois extensions of $\Q $, 
discriminant fields of orthogonally stable characters are not necessarily Galois over $\Q $.
An example is given in \cite[Remark 6.2]{OrthogonalStability} for the absolutely irreducible
characters of degree 56 and 120 of the first Janko group $J_1$. 
However Corollary \ref{abelian} shows that 
all discriminant fields of orthogonally stable characters of finite solvable groups 
are abelian extensions of $\Q $. 

\begin{lemma} \label{squares} 
	Let $K$ be a number field and $\delta \in K^{\times} $. 
	Put $L=K[\sqrt{\delta}]$. 
	Then $(L^{\times})^2  \cap K^{\times} = (K^{\times})^2 \cup  \delta (K^{\times})^2 $.
\end{lemma}

\begin{proof}
	Any $x\in L^{\times}$ is of the form $x=a+b\sqrt{\delta }$ with $a,b \in K$. 
	Then $$x^2 = (a+b\sqrt{\delta })^2 = a^2+b^2\delta + 2ab \sqrt{\delta }  .$$ 
	So $x^2 \in K$ if and only if $2ab=0$ and then $x^2=a^2 \in (K^{\times})^2$ or 
	$x^2=b^2\delta \in \delta (K^{\times})^2$.
\end{proof}

\subsection{Odd square classes}

A field extension $L/K$ is said to be {\em unramified at 2} 
if no prime ideal of $K$ that contains $2$ is ramified in $L/K$. 
\\
If $F\subseteq K \subseteq L \subseteq M$ is a tower of number fields and 
$M/F$ is unramified at 2 then clearly also $L/K$ is unramified at 2
(see for instance \cite[Kapitel III \S 2]{Neukirch}).

Recall that a square class $\delta (K^{\times})^2$ is called {\em odd} 
if $\nu (\delta (K^{\times})^2) $ is even for all dyadic valuations $\nu : K^{\times}\to \Z $ (see 
Definition \ref{oddsc}) 

\begin{lemma} \label{unramodd}
	\begin{itemize}
		\item[(a)] If $\delta \in (K^{\times})^2 $ then $\delta (K^{\times})^2$ is odd. 
		\item[(b)] If $\delta \in K^{\times} \setminus (K^{\times})^2$ and 
			$2$ is unramified in $K[\sqrt{\delta}]/K$ then 
			$\delta (K^{\times})^2$ is odd.
		\item[(c)] The example $K=\Q $ and $\delta = -1$ shows that 
			the converse of (b) is not true. 
	\end{itemize} 
\end{lemma}

\begin{proof}
	The only statement that might require a proof is (b). 
There is a bijection between prime ideals $\wp $ of $K$ and 
	surjective valuations $\nu = \nu_{\wp } : K^{\times} \to \Z $. 
	The prime ideal $\wp $ is unramified in a field extension $L/K$ if 
	for any extension $\tilde{\nu }$ of  $\nu _{\wp }$ 
	we have $\tilde{\nu }(L^{\times}) = \nu _{\wp} (K^{\times})$ (see for instance \cite[Kapitel II \S 8]{Neukirch}).
	For the particular case that $L=K[\sqrt{\delta}]$ this tells us that 
	$$\nu (\delta ) = 2\tilde{\nu }(\sqrt{\delta}) \in 2 \Z $$
	for all dyadic valuations $\nu :K^{\times}\to \Z$, so $\delta (K^{\times})^2$ is odd.
\end{proof}

\begin{lemma} \label{NormSQ}
	Let $L/K$ be a finite extension of number fields 
	and let $\delta \in L^{\times}$.
	Then $N_{L/K}(\delta (L^{\times})^2) \subseteq N_{L/K}(\delta ) (K^{\times})^2$. 
	So the norm of a square class of $L$ defines a square class of $K$.
	\\
	If $\delta (L^{\times})^2 $ is an odd square class of $L^{\times}$ then also 
	$N_{L/K}(\delta ) (K^{\times})^2$ is  odd.
\end{lemma}

\begin{proof}
	Assume that $\delta (L^{\times})^2$ is an odd square class of $L^{\times}$ and
	let $\nu : K^{\times} \to \Z $ be a surjective dyadic valuation. 
	Let $\tilde{\nu }_i: L^{\times} \to \Q $ ($1\leq i \leq s$) be 
	the distinct extension of $\nu $ and $e_i \in \N $ be
	such that $\nu_i := e_i \tilde{\nu }_i :L^{\times} \to \Z $ is 
	surjective. 
	Then by assumption $\nu _i (\delta (L^{\times})^2 ) = 2 \Z $. 
	By \cite[Kapitel III, Satz (1.2) (iv)]{Neukirch} the 
	valuation of the norm of $\delta $ is a linear combination 
	of $\nu _i(\delta )$, 
	$\nu (N_{L/K}(\delta )) = \sum _{i=1}^s f_i \nu _i(\delta ) $.
	In particular $\nu (N_{L/K}(\delta )) $ is even and hence 
	$N_{L/K}(\delta ) (K^{\times})^2$  is an odd square class of $K$.
\end{proof}

\subsection{Orthogonally simple characters} 

Any orthogonally stable character is the 
sum of orthogonally simple characters (see \cite[Section 4.3]{OrthogonalStability}).
There are three kinds of orthogonally simple characters $\chi $: 
\begin{itemize}
\item[(a)] $\chi = 2 \psi $ for some absolutely irreducible real character $\psi $ with  indicator $-$.
\item[(b)] $\chi = \psi + \overline{\psi} $ for an absolutely irreducible non-real character $\psi  $.
\item[(c)] $\chi $ is an absolutely irreducible real character 
	with  indicator $+$.
\end{itemize}

\begin{proposition} (see \cite[Theorem 4.10]{OrthogonalStability}) \label{ODOS}
	Let $\chi $ be an orthogonally simple character. 
	\begin{itemize}
		\item[(a)] If 
		$\chi = 2 \psi $ for some absolutely irreducible real character $\psi $ with  indicator $-$ 
			then $\chi $ is orthogonally stable with orthogonal discriminant 1. 
		\item[(b)] 
		If $\chi = \psi  + \overline{\psi} $ for an absolutely irreducible non-real character
		$\psi  $ then $\chi $ is orthogonally stable. 
		Write 
		$\Q(\psi ) = \Q(\chi ) [\sqrt{\delta }] $ for $\delta \in \Q(\chi )$. 
			Then $\disc (\chi ) = \delta ^{\psi (1)} (\Q(\chi )^{\times})^2 $.
		\item[(c)] 
 If	$\chi $ is an absolutely irreducible real character
			with  indicator $+$ then $\chi $ is orthogonally stable if and only if $\chi (1)$ is even.
	\end{itemize}
\end{proposition}

\begin{remark} \label{Galoisaction}
Let $\chi $ be an orthogonally stable character with orthogonal discriminant 
	$\delta (\Q(\chi )^{\times})^2$.  Let $\sigma $ be  a Galois automorphism of $\Q(\chi )$. 
Then also $\chi ^{\sigma }$ is an  orthogonally stable character and has 
	orthogonal discriminant $\delta ^{\sigma } (\Q(\chi )^{\times})^2 $. 
\end{remark}

If a sum of orthogonally stable characters $\chi = \chi_1+\ldots + \chi_s $ (say with the
same character field $K = \Q(\chi _i)$ for all $i$) has a smaller character field $\Q(\chi ) \subset K$
then the Galois group $\Gal (K/\Q(\chi ))$ acts on the set $\{ \chi_1,\ldots ,\chi_s\}$.
By Lemma \ref{NormSQ} we hence can choose suitable representatives $\delta _i \in \disc (\chi _i)$ 
of the orthogonal discriminants, such that 
$\disc (\chi ) = \prod_{i=1}^s \delta _i (\Q(\chi )^{\times})^2 $. 
In this sense we get the following result:

\begin{theorem} \label{OGNorm} (see \cite[Theorem 4.13]{OrthogonalStability} for a precise statement)
If $\chi $ is an orthogonally stable character written as the sum 
$\chi = \chi_1+\ldots + \chi_s $ of orthogonally simple characters $\chi _i$ then 
	all $\chi _i$ are orthogonally stable and 
	$\disc (\chi )  = \prod _{i=1}^s \disc (\chi _i )$ 
	is a product of norms of orthogonal discriminants 
	of orthogonally simple characters.
\end{theorem}

Using Lemma \ref{NormSQ} we hence get the following corollary.

\begin{corollary}\label{simplegenuegt}
	In the notation of Theorem \ref{OGNorm} the orthogonal discriminant of $\chi $ is 
	odd if the orthogonal discriminants of all the $\chi _i$ are odd.
\end{corollary}

\section{Proof of Theorem \ref{main}} 

By Corollary \ref{simplegenuegt} it is enough to prove Theorem \ref{main} 
for orthogonally stable orthogonally simple characters. 
There are three cases (a), (b), (c) as given in Proposition \ref{ODOS}.
Note that in cases (a) and (b) the orthogonal discriminants 
are odd for arbitrary finite groups.
For case (a) this is already stated in Proposition \ref{ODOS}:

\begin{lemma}
	If $\chi $ is an orthogonal simple character as in Proposition \ref{ODOS} (a) 
	then $\disc(\chi ) = 1$. In particular 
	Theorem \ref{main} holds in case (a).
\end{lemma} 

For case (b) we prove the following lemma:

\begin{lemma} 
	If $\chi $ is an orthogonal simple character as in Proposition \ref{ODOS} (b) 
	then the orthogonal discriminant of $\chi $ is odd. In particular 
Theorem \ref{main} holds  in case (b). 
\end{lemma}

\begin{proof}
Let $K:= \Q(\chi )$ and $L:= \Q(\psi )$. 
	Then $L= K[\sqrt{\delta }] $ for some totally negative $\delta \in K^{\times}$ and $K$ is the
maximal real subfield of $L$. 
	By Proposition \ref{ODOS} 
	the orthogonal discriminant of $\chi $ is $\delta ^{\psi (1)} (K^{\times})^2 $. 
In particular it is a square if $\psi(1)$ is even. 
\\
	So assume that $\psi (1)$ is odd.
Then  Theorem \ref{odddegcor} states that  either the conductor of $\psi $
	is odd or
$\Q[\sqrt{-1}] \subseteq L$.
	In the latter case  $L=K[\sqrt{-1}]$ and hence $\delta = -1 $. 
	\\
So we are left with the case that the conductor of $\psi $ 
is odd and hence $L/K$ is unramified at 2.
	Then Lemma \ref{unramodd} (b) shows that the orthogonal discriminant of $\chi $ is
	odd.
\end{proof}

\begin{remark} 
	To prove the statement of Theorem \ref{main} for 
	arbitrary groups $G$, it is enough to show that 
	all indicator + even degree absolutely irreducible characters
	have odd orthogonal discriminant. 
\end{remark}

\begin{proof} (of Theorem \ref{main}) 
	We now turn to the proof of Theorem \ref{main}. 
As we have shown the Theorem for case (a) and (b) we assume that 
we are in case (c), so $\chi $ is an absolutely irreducible real character of 
Frobenius Schur indicator +.
As $\chi $ is orthogonally stable the degree of $\chi $ is even.
\\
Now we need to assume that $G$ is a solvable group. 
Let $G$ be a minimal counterexample. 
\\
By the minimality of $G$ 
the restriction of $\chi $ to any proper subgroup of $G$ is
not orthogonally stable.
Now $G$ is solvable, so it has a normal subgroup $N\unlhd G$ of prime index $p:=[G:N]$. 
\\
If $p$ is odd then $\chi _{|N} = \chi_1+\ldots +\chi _p$ for absolutely irreducible real 
characters $\chi _i$ of even degree. 
In particular $\chi _{|N}$ is orthogonally stable, a contradiction. 
\\
So $[G:N] = 2$ and $\chi _{|N } = \chi _1 + \chi _2$. 
As $\chi _{|N} $ is not orthogonally stable we conclude that $\chi _1$ and $\chi _2$ 
are real characters of odd degree and hence by Corollary \ref{realodd} the conductor of 
	$\chi _1$ is odd. 
	\\
	Put $K:=\Q(\chi )$, $L:=\Q(\chi_1) = \Q(\chi _2)$.
	Then 
	$L/K$ is unramified at 2.
	\\
Let $\rho $ be a representation of $N$ affording 
the character $\chi _1$ 
and  $G= \langle N , h \rangle $ with $h^2 \in N$.
Then $\chi = \chi _1 ^{G}  $ is induced from the normal subgroup $N$.
So over $L$ we can write the representation $\Rho $ of $G$ with 
character $\chi $ as 
$$\Rho(g) =  \diag (\rho (g) , \rho (g^h) ) \mbox{ for all }  g\in N \mbox{ and } \Rho (h) =  
\left( \begin{array}{cc} 0 & 1 \\ \rho(h^2) & 0 \end{array} \right) .$$
In particular the $\Rho (G)$-invariant forms are of the form 
$\diag (F,F)$ for the  $\rho (N)$-invariant forms $F$. 
This shows that the orthogonal discriminant of 
$\chi $ is a square in $L$. 
\\
	If $K=L$ then $\disc(\chi ) = (K^{\times})^2$. 
\\
If $[L:K] = 2 $ then $L=K[\sqrt{\delta } ]$ for some totally positive $\delta \in K$.
	As $2$ is unramified in $L$, Lemma \ref{unramodd} (b) shows that $\delta (K^{\times})^2$ is 
	an odd square class of $K$. 
	Then 
	$\disc(\chi) \in (L^{\times})^2 \cap K = (K^{\times})^2 \cup \delta (K^{\times})^2 $ (see Lemma \ref{squares})
	is odd. 
\end{proof}

\begin{corollary}\label{abelian}
	If $\chi $ is an orthogonally simple orthogonally stable character of a
	solvable group $G$ then the discriminant field $\Deltaa (\chi )$ is the 
	character field of some subgroup of $G$. 
	\\
	The discriminant field of an orthogonally stable character of a solvable group
	is a subfield of the compositum of the discriminant fields of orthogonally 
	simple characters and 
	in particular an abelian number field.
\end{corollary}

\section{$p$-groups} \label{pgrps}

This section derives  explicit formulas for the  orthogonal discriminants
of orthogonally stable characters of $p$-groups. 

\subsection{Odd primes, the statement}

Let $p$ be an odd prime and $G$ be a finite $p$-group. 
Let $\zeta = \zeta _p $ be a primitive $p$-th root of unity,
 $Z := \Q [\zeta ]$ the $p$-th cyclotomic field with maximal 
 real subfield $Z^+ = \Q[\zeta + \zeta ^{-1}]$ 
 of index 2 in $Z$.

 \begin{definition} 
	 Let $\delta _p\in Z^+$ be such that 
	 $Z=Z^+[\sqrt{\delta _p}]$.
 \end{definition} 

 If $p\equiv 3 \pmod{4} $ then $\delta _p = -p$ is one possible choice. 
 For arbitrary $p$ we can choose 
 $$
\delta_p  = (\zeta -\zeta^{-1})^2 = -N_{Z/Z^+}(1-\zeta^2) = \zeta^2+\zeta^{-2}-2 \in Z^+,$$
 a totally negative generator of the prime ideal of $Z^+$ that
 divides $p$.

 Let 
$\chi $ be an orthogonally stable character of $G$. 
As the only real character of $G$ is the trivial character, 
this is equivalent to the fact that 
the character field $K:=\Q(\chi )$ is real and $\chi $ does not
contain the trivial character of $G$ as a constituent. 
Clearly the conductor of $\chi $ is a power of $p$, say $p^f$,
whence 
 $K$ is a subfield of the cyclotomic field $Z_f:= \Q [\zeta _{p^f}]$,
 a cyclic extension of degree $p^{f-1}(p-1)$ over $\Q $.
As $f$ is minimal, the index of $K$ in $Z_f$ is prime to $p$ and
hence $[Z_f:K] = a$ 
where $a$ divides $p-1$. 

 \begin{lemma}\label{adivdeg}
$a$ is even and divides $\chi (1)$.
\end{lemma} 

Denote by 
$Z_f^+:= \Q[\zeta_{p^f} + \zeta _{p^f}^{-1} ]$ the maximal real subfield of 
$Z_f$.
Then also $Z_f = Z_f^+[\sqrt{\delta_p }] $.
As $K$ is real, we have $K \subseteq Z_f^+$, in particular $a$ is even.
Put 
$$\delta _K  := N_{Z_f^+/K}(\delta_p) \in K \cap Z^+.$$

 \begin{theorem}\label{mainP}
The orthogonal discriminant of $\chi $ is $\delta_K ^{\chi (1)/a} (K^{\times })^2$. 
\end{theorem}

\begin{corollary}
         If $p\equiv 3 \pmod{4} $ then the orthogonal discriminant of
	 $\chi $ is $(-p)^{\chi (1)/2 }(K^{\times})^2 $.
 \end{corollary}

 \subsection{Odd primes, a proof} 

 As before let $Z_f:= \Q[\zeta _{p^f}]$ be the $p^f$-th cyclotomic field.
 Then $Z_f$ is a cyclic Galois extension of $\Q $ with 
 Galois group 
 $$ \Gamma _f \cong (\Z/p^f \Z )^{\times } \cong C_{p-1} \times C_{p^{f-1}} .$$ 
 In particular $\Gamma _f$ contains a unique subgroup 
 $\Gamma _1 =\langle \gamma _f \rangle $ of order $p-1$.  

 \begin{remark}
	 For $f_1\leq f_2$  the restriction 
	 of $\gamma _{f_2}$ to $Z_{f_1} $ generates the group $\langle \gamma _{f_1} \rangle $.
 \end{remark}

 Write $\chi = \sum_{i=1}^s n_i \psi _i  $ with $n_i \in \N$ and $\psi _i\in 
 \Irr (G)$ for
 pairwise distinct complex irreducible characters of $G$. 
 As $\chi $ is orthogonally stable, none of the $\psi _i$ is the trivial character and hence
 the character field $\Q (\psi _i ) = Z_{f_i}$ for suitable $f_i\in \N$ by Theorem \ref{pgrpchar}.
 Let $F $ be the maximum of the $f_i$. 
 As the character field $K = \Q(\chi )$ has index $a $ dividing $p-1$ in $Z_f$ 
 we have $f\leq F$. Put $\gamma := \gamma _F^{(p-1)/a}$, so that the 
 restriction of $\gamma $ to $Z_{f_i}$ generates the subgroup of order 
 $a$ in $\Gamma _{f_i}$. 
 In particular
 $\langle \gamma \rangle $ acts on the set of $\psi _i$ with orbits 
 of length $a$. 
 Assume that $\psi _1,\ldots , \psi _t$ represent these different orbits and put
$\chi _i := (\sum _{j=1}^a \psi_i ^{\gamma ^j} )$ for $1\leq i \leq t$.
Then
$$\chi = \sum _{i=1}^t n_i (\sum _{j=1}^a \psi_i ^{\gamma ^j} ) = 
\sum _{i=1}^t n_i \chi _i $$  
where $a\psi _i (1) =  \chi _i(1) $, $\chi (1) = a \sum _{i=1}^t n_i \psi _i(1) $, which proves Lemma \ref{adivdeg}.

\begin{lemma}
	Let $1\leq i \leq t$.
	\begin{itemize}
		\item[(a)] $\Q(\chi _i)$ is the subfield of index $a$ in 
			$\Q(\psi _i)$.
	\item[(b)]
$\Q(\chi _i)/K\cap \Q(\chi _i)$ is a power of $p$.
	\item[(c)] $\delta _K \in \Q(\chi _i)$.
	\item[(d)]
 $\chi_i$ is orthogonally stable of orthogonal discriminant
$\delta _K ^{\psi_i(1)} (\Q(\chi _i)^{\times })^2 $. 
	\end{itemize}
\end{lemma}

 \begin{proof} 
	 (a) By construction $\Q(\chi _i)$ is the subfield of 
	 index $a$ in $\Q(\psi _i) = Z_{f_i} $ and
	 $K$ is the subfield of index $a$ in $Z_f$. 
	 \\
	 (b) If
 $f_i \leq f$ then $\Q(\chi _i) \subseteq K$ and hence
	 $\Q(\chi _i ) \cap K = \Q(\chi _i)$.
	 If $f_i \geq f$ then $K\subseteq \Q(\chi _i) $ of 
	 index $p^{f_i-f}$. 
	 \\
 To see (c) we remark that 
 $\delta _K  = N_{Z_f^+/K}(\delta_p) \in K \cap Z^+ $ and 
	 $K\cap Z^+$ is the subfield of index $a$ in $Z=Z_1$.
As $Z$ is a subfield of all the character fields $\Q(\psi _i)$ 
also its subfield $K\cap Z^+$  of index $a$ is contained in the subfield 
	 $\Q(\chi _i)$.
 \\
	 For part (d) we use the product formula from Theorem \ref{OGNorm}.
        We first remark that
        for $1\leq i \leq t$ the orthogonally simple character
        $\psi_i + \overline{\psi_i}$
has  orthogonal discriminant
        $\delta = \delta_p^{\psi_i(1)}  \in Z^+ \leq Z_{f_i}^+$.
        By
        \cite[Proposition 4.11]{OrthogonalStability}  the
        orthogonal discriminant of $\chi _i$ is obtained by
        taking the norm
	 of $\delta $ as $\disc(\chi_i) = \delta _K^{\psi_i(1)} \in Z^+\cap K\subseteq \Q(\chi _i)$. 
\end{proof}

Now $\Q(\chi _i)/K\cap \Q(\chi _i) $ is odd and 
$\delta _K \in K$,
so we now can use Theorem \ref{OGNorm}
 to compute
$$\disc (\chi ) = \prod _{j=1}^t (\delta _K ^{n_i \psi_i(1)}) (K^{\times })^2 =
\delta _K^{\sum_{i=1}^t n_i \psi_i(1) } (K^{\times })^2 = \delta_K^{\chi(1)/a} 
(K^{\times})^2 $$ 
	which proves Theorem \ref{mainP}.

\subsection{$2$-groups}

	\begin{theorem}\label{main2}
        Let $G$ be a $2$-group and $\chi $ be an orthogonally stable
        orthogonal character of $G$.
        Then the degree $\chi (1)$ is even and
        the orthogonal discriminant of $\chi $ is
        $\disc(\chi ) = (-1)^{\chi(1)/2} $.
        \end{theorem}

	\begin{proof}
                As the discriminant of an orthogonally stable character
                is the product of the discriminants of all orthogonally
		simple summands by Theorem \ref{OGNorm}, 
                it is clearly enough to show the statement for orthogonally
                simple characters $\chi $.
                We proceed by induction on the group order.
                If $|G| = 2$ then there are no orthogonally stable
                characters, so here the statement is trivial.
                If $|G| = 4$ then either $G\cong C_2\times C_2$
                and there are no orthogonally stable characters or
                $G \cong C_4$ with a unique orthogonally simple orthogonally stable character
                $\chi =\psi + \overline{\psi }$ for the complex
                character $\psi $ of degree 1 and with character field
                $\Q(\psi ) = \Q [i] = \Q[\sqrt{-1}] $. Then by
		Proposition \ref{ODOS} (b)  we get
                $\disc(\chi ) = (-1)^{\psi(1)} = (-1) ^{\chi(1)/2 }$.
                \\
   Now let $|G| > 4$ and let $\chi $ be an orthogonally stable
                orthogonally simple character of $G$.
                Without loss of generality we assume that the
                corresponding representation $\rho _{\chi }$ affording
                the character $\chi $ is faithful.
                Choose a normal subgroup $N\unlhd G$ with $|G/N| = 2$.
                If the restriction of $\chi $ to $N$ is orthogonally stable,
                then the statement follows by induction.
                So assume that $\chi _{|N }$ is not orthogonally stable,
                i.e. there is an absolutely irreducible
                real constituent $\psi $ of $\chi _{|N }$ of
                odd degree. But $N$ is a 2-group so $\psi (1) = 1$, and
                $\chi $, being orthogonally simple, is $\psi ^G$ and
                of degree $\chi(1) =2$. As $\psi $ is a real linear
                character the image $\psi (N)$ has order 1 or 2 and
                 $G$ has order 8 and a normal subgroup
		 $N\cong C_2\times C_2$. We conclude that
		 $G \cong D_8$ and $\disc(\chi ) = -1 $.
        \end{proof}

\end{document}